\documentclass[aps,prl,twocolumn,groupedadress]{revtex4-1}

\usepackage{graphicx}
\usepackage{dcolumn}
\usepackage{bm}
\usepackage{dsfont,amsmath,amssymb,color,units,pstricks}

\newcommand \be  {\begin{equation}}
\newcommand \bea {\begin{eqnarray} \nonumber }
\newcommand \ee  {\end{equation}}
\newcommand \eea {\end{eqnarray}}

\def \P{ \mathbb  P }

\def \C{ \mathbb  C }
\def \N{ \mathbb  N }

\begin{document}

\title{Invariant $\beta$-ensembles and the Gauss-Wigner crossover \\}

\author{Romain Allez$^{1,2}$ and Jean-Philippe Bouchaud$^2$ and Alice Guionnet$^3$}
\affiliation{$^1$  Universit{\'e} Paris-Dauphine, Ceremade, 75\,016 Paris, France, }
\affiliation{$^2$ Capital~Fund~Management, 6--8 boulevard Haussmann, 75\,009 Paris, France}
\affiliation{$^3$ U.M.P.A.
ENS de Lyon 
46, all\'ee d'Italie 
69364 Lyon Cedex 07 - France. }
\date{\today}

\begin{abstract}
We define a new diffusive matrix model converging towards the $\beta$ -Dyson Brownian motion for all $\beta\in [0,2]$
that provides an explicit construction of  $\beta$-ensembles of random matrices that is invariant under the orthogonal/unitary group. 
For small values of $\beta$, our process allows one to interpolate smoothly between the Gaussian distribution and 
the Wigner semi-circle.  The interpolating limit distributions form a one parameter family that can be explicitly computed. This also allows us to compute the finite-size corrections to the semi-circle.
\end{abstract}

\maketitle


\vspace{1cm}
\footnotesize

\normalsize

Since Wigner's initial intuition that the statistical properties of the eigenvalues of random matrices should provide a good description of the excited 
states of complex nuclei, Random Matrix Theory has become one of the prominent field of research, at the boundary between atomic physics, solid state physics, 
statistical mechanics, statistics, probability theory and number theory \cite{Handbook,Bai,AGZ}.
It is well known that the joint distribution of the eigenvalues of a large Gaussian random matrix can be expressed as the Boltzmann-Gibbs equilibrium
weight of a one-dimensional repulsive Coulomb gas confined in an harmonic well. However, the effective ``inverse temperature'' $\beta$ of the 
system cannot take arbitrary values but is {\it quantized} (in units of the repulsive Coulomb potential). Depending on the symmetry of the random 
matrix, only three values are allowed $\beta=1$ for symmetric real matrices, $\beta=2$ for Hermitian matrices and $\beta=4$ for the symplectic 
ensemble. This is known as Dyson's ``threefold way''. The existence of matrix ensembles that would lead to other, possibly continuous, values of $\beta$, is
a very natural question, and the quest for such ensembles probably goes back to Dyson himself. Ten years ago, Dumitriu and Edelman \cite{Dumitriu} 
have proposed an explicit construction of tri-diagonal matrices with non-identically distributed elements whose joint law of the eigenvalues 
is the one of $\beta$-ensembles for general $\beta$. Another construction is proposed in \cite[p. 426-427]{Handbook} (see also \cite{Forrester3}) and uses 
a bordering procedure to construct recursively a sequence of matrices with eigenvalues distributed as $\beta$-ensembles. This construction 
gives not just the eigenvalue probability density of one matrix of the sequence but also the joint 
eigenvalue probability density of all matrices. 
This has lead to a renewed 
interest for those ensembles, that have connections with many problems, both in physics and in mathematics, see e.g. \cite{Forrester,AGZ}. 
The aim of the paper is to provide another construction of $\beta$-ensembles that is, at least to our eyes, natural and transparent, 
and respects by construction the orthogonal/unitary symmetry \cite{Vivo}. Another motivation for our work comes from the recent development of free probability theory. 
``Freeness'' for random matrices is the natural extension
of independence for classical random variables. Very intuitively, two real symmetric matrices ${\bf A,B}$ are mutually free in the large $N$ limit if the eigenbasis of ${\bf B}$
can be thought of as a random rotation of the eigenbasis of ${\bf A}$ (see e.g. \cite{Verdu} for an accessible introduction to freeness and for more rigorous statements). 
``Free convolution'' then 
allows one to compute the eigenvalue distribution of the sum ${\bf A+B}$ from the eigenvalue distribution of ${\bf A}$ and ${\bf B}$, much in the same way 
as convolution allows one to compute the distribution of the sum of two independent random variables. In this context, the Wigner semi-circle distribution
appears as the limiting distribution for the sum of a large number of free random matrices, exactly as the Gaussian is the limiting distribution for the 
sum of a large number of {\it iid} (independent and identically distributed) random variables. A natural question, from this perspective, is whether one can build a natural framework that 
interpolates between these two limits. 

Let us first recall Dyson's Brownian motion construction of the GOE \cite{Dyson} (for the sake of simplicity, we will only consider here extensions of the $\beta=1$
ensemble, but similar considerations hold for $\beta=2$ Hermitian matrices   
see \cite{paper_math} for full details). 
One introduces a fictitious time $t$ for the evolution of an $N \times N$ real symmetric matrix ${\bf M}(t)$. The evolution of the symmetric matrix is governed by the following
stochastic differential equation (SDE): 
\begin{equation}
{\rm d} {\bf M}(t) = -\frac{1}{2} {\bf M}(t) {\rm d}t + {\rm d}{\bf H}(t)
\end{equation} 
where $ {\rm d}{\bf H}(t)$ is a symmetric Brownian increment (i.e. a symmetric matrix whose entries above the diagonal are independent Brownian increments with variance 
$\langle {\rm d}{\bf H}_{ij}^2(t) \rangle = \frac{\sigma^2}{2} (1+\delta_{ij}) {\rm d}t$). Standard second order perturbation theory
allows one to write the evolution equation for the eigenvalues $\lambda_i$ of the matrix ${\bf M}(t)$:
\begin{equation}\label{dyson}
{\rm d}\lambda_i = -\frac{1}{2} \lambda_i  {\rm d}t  + \frac{\sigma^2}{2}  \sum_{j \neq i} \frac{{\rm d}t }{\lambda_i - \lambda_j} + \sigma
{\rm d}b_i,
\end{equation}
where $b_i(t)$ are independent standard Brownian motions. 
This defines Dyson's Coulomb gas model, i.e. ``charged'' particles on a line, with positions
$\lambda_i$, interacting via a logarithmic potential, subject to some thermal noise and confined by a harmonic potential. One can deduce from the 
above equation the Fokker-Planck equation for the joint density $P(\{\lambda_i\},t)$, for which the stationary joint probability density function (pdf) is readily found to
be:
\begin{equation}\label{Pstar}
P^*(\{\lambda_i\}) = Z \prod_{i < j} |\lambda_i - \lambda_j|^\beta \exp \left[ - \frac{1}{2 \sigma^2} \sum_i \lambda_i^2 \right],
\end{equation}
with $\beta \equiv 1$, independently of $\sigma^2$ and where $Z$ is a normalization factor. The above expression is the well known joint distribution of the eigenvalues of an $N \times N$ 
random GOE matrix. The Wigner distribution can be recovered either by a careful analysis of the mean marginal univariate distribution $\rho(\lambda)=
 \int \dots \int {\rm d}\lambda_2 \dots {\rm d}\lambda_N P^*(\lambda=\lambda_1,\lambda_2, \dots, \lambda_N)$ in the large $N$ limit \cite{Mehta}, or by using the 
above SDE \eqref{dyson} to derive a dynamical equation for the Stieltjes transform $G(z,t)$ of $\rho(\lambda,t)$:
\be
G(z,t) = \frac{1}{N} \sum_{i=1}^N \frac{1}{\lambda_i(t) - z} , \quad z \in \C.
\ee
With this scaling, the spectrum is spread out in a region of width of order $\sqrt{N}$ and therefore $z\sim \sqrt{N}$ and $G \sim 1/\sqrt{N}$.
Applying It\^o's formula to $G(z,t)$ and using \eqref{dyson}, we obtain the following Burgers equation for $G$\cite{Rogers_Shi}:
\begin{align}
2\frac{\partial G}{\partial t}= \frac{\alpha \sigma^2 N}{2}  \frac{\partial G^2}{\partial z} + 
\frac{\partial z G}{\partial z} + (2 - \alpha) \frac{\sigma^2}{2} \frac{\partial^2 G}{\partial z^2} \label{RS}
\end{align}
where $\alpha$ is introduced for later convenience, with $\alpha=1$ for now. Note that we have neglected a noise term of mean zero and variance of 
order $1/N\sqrt{N}$: $G$ converges to the solution of Eq. (\ref{RS}) only up to fluctuations of order 
$\langle G^2 \rangle-\langle G\rangle^2 \sim 1/N^3$ (in agreement with e.g. \cite[Theorem 9.2]{Bai}). The contribution of these fluctuations 
is thus $1/N$ smaller than the diffusion term in Eq. (\ref{RS}), and can be neglected. 

For large $N$, the last (diffusion) term of Eq. (\ref{RS}) is of order $1/N$ smaller than the other ones. 
To leading order, the stationary solution (where the derivative with respect to time is set to $0$) can be easily integrated with respect to $z$ as: 
\begin{equation}\label{RS_statio}
\frac{1}{2} \alpha \sigma^2 N G^2_\infty(z) + z G_\infty(z) = -1 \,,
\end{equation}
where the integration constant comes from the boundary condition $G(z)\sim -1/z$ when $z\rightarrow \infty$.
It is then easy to solve this equation to find the Stieltjes transform that indeed corresponds to the Wigner semi-circle density:
\begin{align}
&G_\infty(z) = \frac{1}{\alpha\sigma^2N} \left[\sqrt{z^2 - 2\alpha\sigma^2N} - z\right] \notag \\ 
& \longrightarrow \rho(\lambda) = \frac{1}{\pi \alpha \sigma^2N} \sqrt{2\alpha \sigma^2 N- \lambda^2} 
\mathds{1}_{\{|\lambda|\leq \sqrt{2\alpha N} \sigma\}}. \label{Wigner}
\end{align}

Now let us turn to the central idea of the present paper. In Dyson's construction, the extra Gaussian slice ${\rm d}{\bf M}(t)$ that is added to 
${\bf H}(t)$ is chosen to be independent of ${\bf M}(t)$ itself. The eigenbasis of ${\rm d}{\bf H}(t)$ is a random rotation, taken uniformly over the
orthogonal group. As mentioned above, this corresponds to free addition of matrices, and Eq. (\ref{RS}) can indeed be derived (for $N = \infty$) using
free convolution \cite{Verdu}. If instead we choose to add a random matrix ${\rm d}{\bf Y}(t)$ that is {\it always diagonal in the same basis} as that of ${\bf M}(t)$,
the process becomes trivial. The diagonal elements of ${\bf M}(t)$ are all sums of {\it iid} random variables, and the eigenvalue distribution converges towards 
the Gaussian. The construction we propose is to alternate randomly the addition of a ``free'' slice and of a ``commuting'' slice. 
More precisely, our model is defined as follows: we divide time into small intervals of length $1/n$ and for each interval $[k/n;(k+1)/n]$, we choose independently 
Bernoulli random variables $\epsilon_k^n, k\in\N$ such that $\P[\epsilon_k^n=1] = p = 1-\P[\epsilon_k^n=0]$. Then, setting $\epsilon_t^n = \epsilon_{[nt]}^n$, our diffusive matrix 
process simply evolves as: 
\begin{equation}\label{def_model}
{\rm d}{\bf M}_n(t) = -\frac{1}{2} {\bf M}_n(t) {\rm d}t + \epsilon_t^n {\rm d}{\bf H}(t) + (1-\epsilon_t^n) {\rm d}{\bf Y}(t)
\end{equation} 
where ${\rm d} {\bf H}(t)$ is a symmetric Brownian increment as above and where ${\rm d}{\bf Y}(t)$ is a symmetric matrix that is co-diagonalizable with ${\bf M}_n(t)$ 
(i.e. the two matrix have the same eigenvectors)
but with a spectrum given by $N$ independent Brownian increments of variance ${\sigma^2}{\rm d}t$.
It is clear that the eigenvalues of the matrix ${\bf M}_n(t)$ will cross at some points but only in intervals $[k/n;(k+1)/n]$ for which $\epsilon_k^n=0$ (in the other intervals 
where they 
follow Dyson Brownian motion with parameter $\beta=1$, it is well known that the repulsion is too strong and that collisions are avoided). In such a case, the eigenvalues are re-numbered 
at time $t=(k+1)/n$ in increasing order. 


Now, using again standard perturbation theory, it is easy to derive the evolution of the eigenvalues of ${\bf M}_n(t)$ denoted as $\lambda_1^n(t)\leq \dots \leq \lambda_N^n(t)$:
\begin{equation}\label{dyson_epsilon}
{\rm d}\lambda_i^n = -\frac{1}{2} \lambda_i^n  {\rm d}t  + \epsilon_t^n \frac{\sigma^2}{2}  \sum_{j \neq i} \frac{{\rm d}t}{\lambda_i^n - \lambda_j^n} + \sigma 
{\rm d}b_i
\end{equation}
where the $b_i$ are independent Brownian motions also independent of the $\epsilon_k^n, k \in \N$.

A mathematically rigorous derivation provided in
\cite{paper_math} allows one to show that the scaling limits $\lambda_i(t)$, when $n \rightarrow \infty$, of the eigenvalues $\lambda_i^n(t)$ obey the 
following modified Dyson SDE:
\begin{equation}\label{dyson_p}
{\rm d}\lambda_i = -\frac{1}{2} \lambda_i  {\rm d}t  + p \frac{ \sigma^2}{2}  \sum_{j \neq i} \frac{{\rm d}t}{\lambda_i - \lambda_j} + \sigma 
{\rm d}b_i,
\end{equation}
with the additional ordering constraint $\lambda_1(t)\leq \dots\leq\lambda_N(t)$ for all $t$. 
One of the difficulty of the proof comes from the fact that when $p < 1$, there is a positive probability for eigenvalues to collide in finite time (the ordering constraint is therefore 
useful at those points to re-start). The idea is then 
to show that collisions are in a sense sufficiently rare for the above SDE to make sense (see \cite{paper_math,Cepa} for further details). Using the 
SDE $\eqref{dyson_p}$, one can derive as above the stationary distribution for the joint distribution of eigenvalues, which is still given by Eq. (\ref{Pstar}) but 
with now $\beta=\alpha=p \leq 1$. A very similar construction can be achieved in the GUE case, leading to $\beta = 2p$. As announced, our dynamical 
procedure, that alternates standard and free addition of random matrices, can lead to any $\beta$-ensemble with $\beta \leq 2$. 
The corresponding matrices ${\bf M}(t)$ are furthermore {\it invariant} under the orthogonal (or unitary) group. This is intuitively clear, since 
both alternatives (adding a free slice or adding a commuting slice) respect this invariance, and lead to a Haar probability measure for the 
eigenvectors (i.e. uniform over the orthogonal/unitary group). We have also proved that a collision leads to a complete randomization of the eigenvectors within the
two-dimensional subspace corresponding to the colliding eigenvalues, see again \cite{paper_math}.

%

\begin{figure}[h!btp] 
	\center
		\includegraphics[scale=0.3]{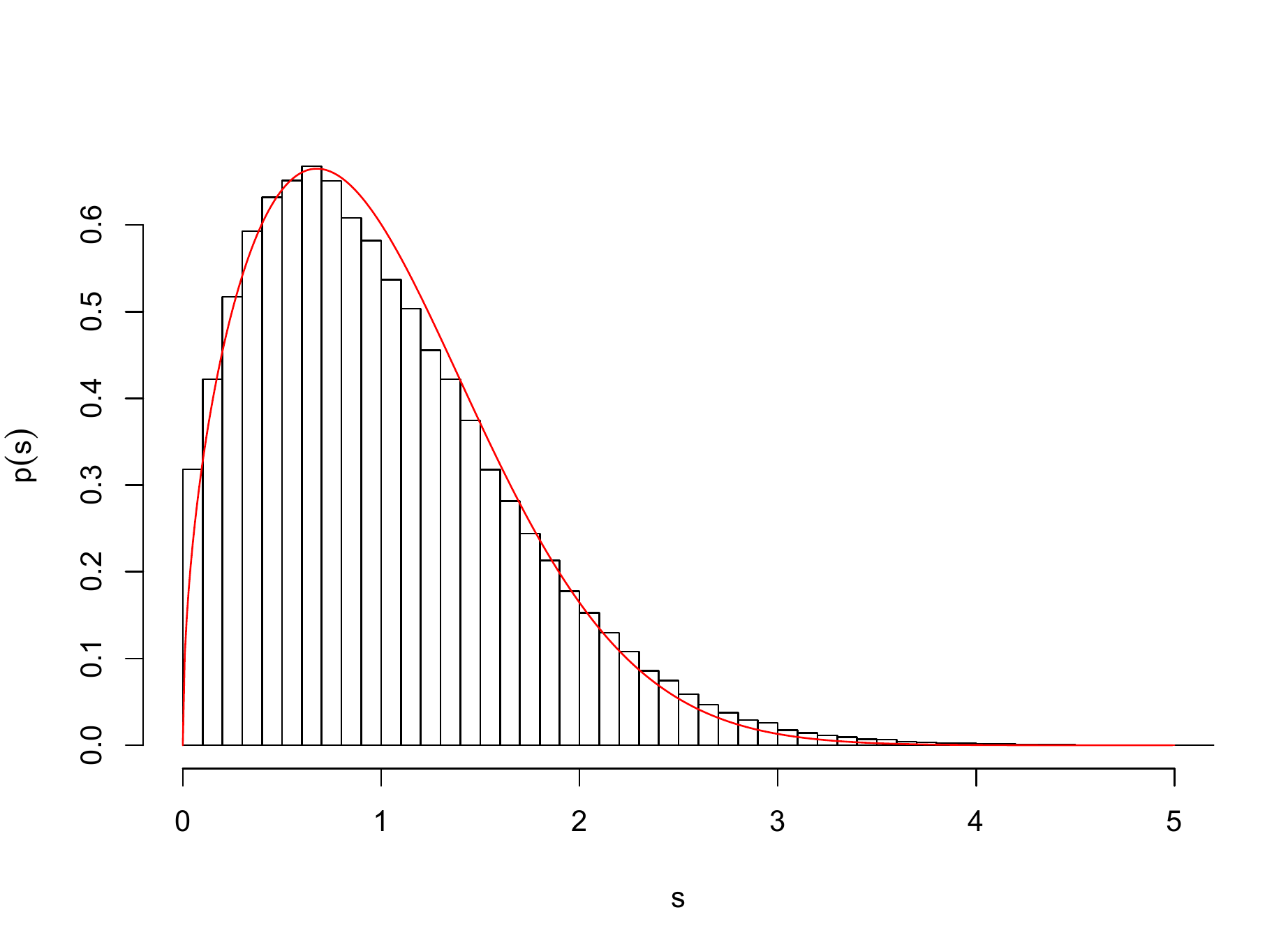}
		\caption{Empirical NNSD $P(s)$ for the matrix ${\bf M}_n(t=\infty)$ for $\beta=p=1/2$ with 
		the Wigner surmise (red curve) corresponding to $\beta=\frac{1}{2}$, which behaves as $s^\beta$ when $s \to 0$. }\label{simul_p_undemi}
\end{figure}

It is well known that the eigenvalue density corresponding to the measure $P^*$ given by \eqref{Pstar} is the Wigner semi-circle for any $\beta > 0$. 
In fact, using \eqref{RS} with now $\alpha=\beta=p$, one immediately finds that the eigenvalue density is a semi-circle with edges at $\pm \sigma \sqrt{2 \beta N}$.
We simulated numerically the matrix ${\bf M}_n(t)$ with $N=200$ for a very small step $1/n$ and until a large value of $t$ so as to reach the stationary distribution 
for the eigenvalues. Then we started recording the spectrum and the nearest neighbor spacings (NNS) every $100$ steps so as to sample the ensemble. We verified that the 
spectral density of ${\bf M}_n(t=\infty)$ is indeed in very good agreement with the Wigner semi-circle distribution for $\beta=1/2$. 
Our sample histogram for the NNS distribution is displayed in Fig.\ref{simul_p_undemi}. 
We also added  the corresponding 
Wigner surmise (which is expected to provide a good approximate description of the NNSD).


From the point of view of a cross-over between the standard Gaussian central limit theorem for random variables and the Wigner central limit theorem for
random matrices, we see that as soon as the probability $p$ for a non-commuting slice is positive, the asymptotic density is the Wigner semi-circle, with a width 
of order $\sqrt{p N}$. A continuous cross-over therefore takes place for $p=2c/N$ with $c$ strictly positive and independent of $N$.  When $c=0$, 
$\rho(\lambda)$ is a Gaussian of rms $\sigma$, which indeed corresponds to the solution of Eq. (\ref{RS}) for $\alpha=0$. 
Setting for simplicity $\sigma=1$, the SDE for the system $({\lambda}_i(t))$ becomes 
\begin{equation}\label{dyson_p_hat}
{\rm d}{\lambda}_i = -\frac{1}{2} {\lambda}_i  {\rm d}t  + \frac{c}{N}  \sum_{j \neq i} \frac{{\rm d}t}{{\lambda}_i - {\lambda}_j} +
{\rm d}b_i,
\end{equation}
with the additional ordering constraint ${\lambda}_1(t)\leq \dots\leq{\lambda}_N(t)$ and the stationary joint pdf is still given by \eqref{Pstar} 
but with now a vanishing repulsion coefficient $\beta = 2c/N$. In
order to elicit the cross-over, we study Eq. \eqref{RS} with $\alpha=2 c/N$.
The stationary differential equation corresponding to \eqref{RS} (note this time that all terms are of the same and the second derivative term is not negligible) 
can be integrated 
with respect to $z$ again as: 
\begin{equation}\label{eq_G_hat}
c  G^2 +  z  G + \frac{{\rm d}  G}{{\rm d}  z} = - 1,
\end{equation}
where the integration constant comes from the boundary condition $ G \sim - 1/ z$ for $ z \to \infty$. Note that 
\eqref{eq_G_hat} can be recovered directly from the saddle point equation route: under the measure ${P}^*$ with $\beta=2c/N$, the energy of a configuration of the $\lambda_i$'s 
can be expressed in term of the continuous state density $\rho$, neglecting terms $\ll 1$, as:
\begin{align*}
\mathcal{E}[\rho]= \frac{1}{2} \int \lambda^2 \rho(\lambda) {\rm d}\lambda 
- c \int\int \ln(|\lambda-\lambda'|) \rho(\lambda) \rho(\lambda') {\rm d}\lambda {\rm d}\lambda'\,.
\end{align*}
The probability density ${P}^*$ therefore rewrites in term of $\rho$ as:
\begin{align*}
{P}^*[\rho] = Z \exp\left(-N \left[\mathcal{E}[\rho]+\int \rho \ln(\rho)  \right] \right) \delta(\int \rho - 1)\,,
\end{align*}
where the entropy term, which is negligible when $\beta=p$ is of order $1$, is now of the same order as the energy term (see \cite{Satya-Dean} for a detailed discussion on the
origin of the entropy term). We now need to minimize 
the quantity $\mathcal{E}[\rho]+\int \rho \ln(\rho)$ with respect to $\rho$. It is easy to see that the unique minimizer $\rho_c$ satisfies: 
\begin{align*}
\int \frac{\lambda \rho_c(\lambda)}{\lambda-z}{\rm d}\lambda - 2c \int\int\frac{\rho_c(\lambda)\rho_c(\lambda')}{(\lambda-z)(\lambda-\lambda')}{\rm d}\lambda {\rm d}\lambda' \\
+ \int \frac{\rho_c'(\lambda)}{\lambda-z}{\rm d}\lambda + \nu = 0
\end{align*} 
where $\nu$ is an integration constant. It is now straightforward to derive \eqref{eq_G_hat} from this last equation by identifying each term and choosing the constant $\nu$ so as to have 
the correct boundary condition for the Stieltjes transform of a probability measure. As expected physically, the diffusion term in \eqref{eq_G_hat} corresponds exactly to the entropy 
contribution to the saddle-point. 

Now, to solve \eqref{eq_G_hat}, we first set ${G}( z):=u'(z)/cu( z)$ to obtain a second order equation on $u$: 
\begin{equation}\label{eq_u}
u''({z})+{z} u'({z}) + c u({z}) = 0\,.
\end{equation}
It follows from the asymptotic behavior of ${G}( z)$ that, for $|{ z}| \rightarrow \infty$, 
\begin{equation}\label{asympt_u}
u({z}) \sim \frac{A_1}{{z}^{c}}\,.
\end{equation}
Eq. \eqref{eq_u} can in turn be transformed 
with the change of function $u( z):=e^{-{z}^2/4} y({z})$ into 
a Schrodinger equation (sometimes called {\it Weber differential equation}) on $y({z})$: 
\begin{equation}\label{eq_y}
y''({z}) + [c-\frac{1}{2} - \frac{1}{4} {z}^2 ] y({z}) = 0\,. 
\end{equation} 
The solutions of \eqref{eq_y} are known (see \cite[page 1067]{handbook_math_funct}) to write as $y({z}) = A_2 D_{c-1}({z}) + A_3 D_{-c}(i{z})$ where $D_{c-1}, D_{-c}$ are parabolic cylinder functions 
and where $A_2$ and $A_3$ 
are two constants. The general solution for $u$ therefore is $u( z)= e^{-{z}^2/4} (A_2 D_{c-1}({z}) + A_3 D_{-c}(i{z}))$ and the correct asymptotic behavior of $u$ is fulfilled for $A_2=0$. 
Now, one can recover the spectral density $\rho_c({\lambda})$ associated to $ G$ by the classical inversion formula: 
\begin{align}
\rho_c({\lambda}) = \lim_{\epsilon\rightarrow 0} \frac{1}{\pi} \Im({G}({\lambda}-i\epsilon))\\
=\frac{1}{c\pi} \frac{1}{|y({\lambda})|^2} (y'_2 y_1 - y_2 y'_1)({\lambda})\,,
\end{align}
where $y_1$ and $y_2$ are respectively the real and imaginary part of $y$. Using \eqref{eq_y}, it is straightforward to see that the derivative of $(y'_2 y_1 - y_2 y'_1)$ 
is identically zero and therefore the function is constant. This gives up to a normalizing constant the exact expression for $\rho_c(\lambda)$. It happens that this constant can be 
computed explicitly \cite{satya}. Eq. \eqref{eq_G_hat} was 
also studied in detail by Kerov \cite{Kerov} and Askey \& Wimp \cite{AskeyWimp} (see also \cite{GaussWigner}). 
The final result for $\rho( \lambda)$ reads, for all $c>0$:  
\begin{align}
\rho_c( \lambda)&= \frac{1}{\sqrt{2\pi} \Gamma(1+c)} \frac{1}{|D_{-c}(i  \lambda)|^2}; \\ 
D_{-c}(z)& = \frac{e^{-z^2/4}}{\Gamma(c)} \int_0^\infty {\rm d}x 
e^{-zx -\frac{x^2}{2}} x^{c-1}.
\end{align}
This expression was again checked with numerical simulations with very good agreement.  
The integral representation for $D_{-c}(z)$ does not hold for $c=0$, but the function $D_{-c}(iu)$ is still well defined for all $c\in (-1;0]$ (see \cite[Theorem 8.2.2]{AskeyWimp}). 
It is easy to check that $\rho_0(u)=e^{-u^2/2}/\sqrt{2\pi}$ when $c=0$, as expected.
When $c\rightarrow \infty$, the Wigner semi-circle law is recovered 
\begin{equation}
\rho_c(u) \approx \frac{1}{2\pi c} \sqrt{4c-u^2}\,.
\end{equation}
At least when $c$ is a non-negative integer, the integral form of $D_{-c}(i\lambda)$ can be computed analytically. This enables to 
find the tails of $\rho_c(u)$ for some values of $c$ and large $u$. The asymptotic behavior reads:
\begin{equation}
\rho_c(u) \sim u^{2c} e^{-u^2/2}\quad (|u|\to \infty). 
\end{equation}

\begin{figure}[h!btp] 
	\center
		\includegraphics[scale=0.4]{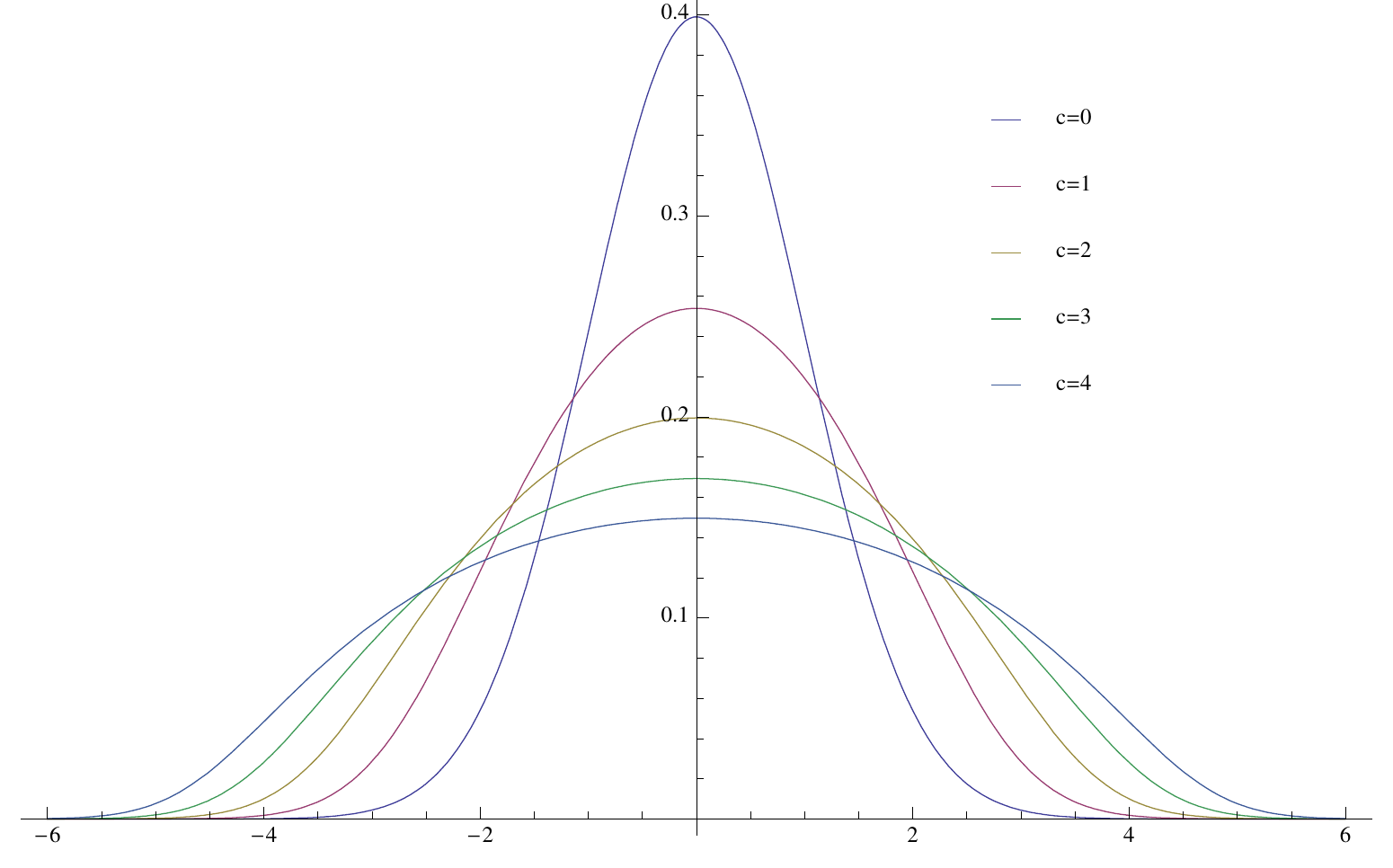}
		\caption{Density $\rho_c(u)$ for $c=0,1,2,3,4$, showing the progressive deformation of the Gaussian towards Wigner's semi-circle.}
\end{figure}

\begin{figure}[h!btp] 
	\center
		\includegraphics[scale=0.35]{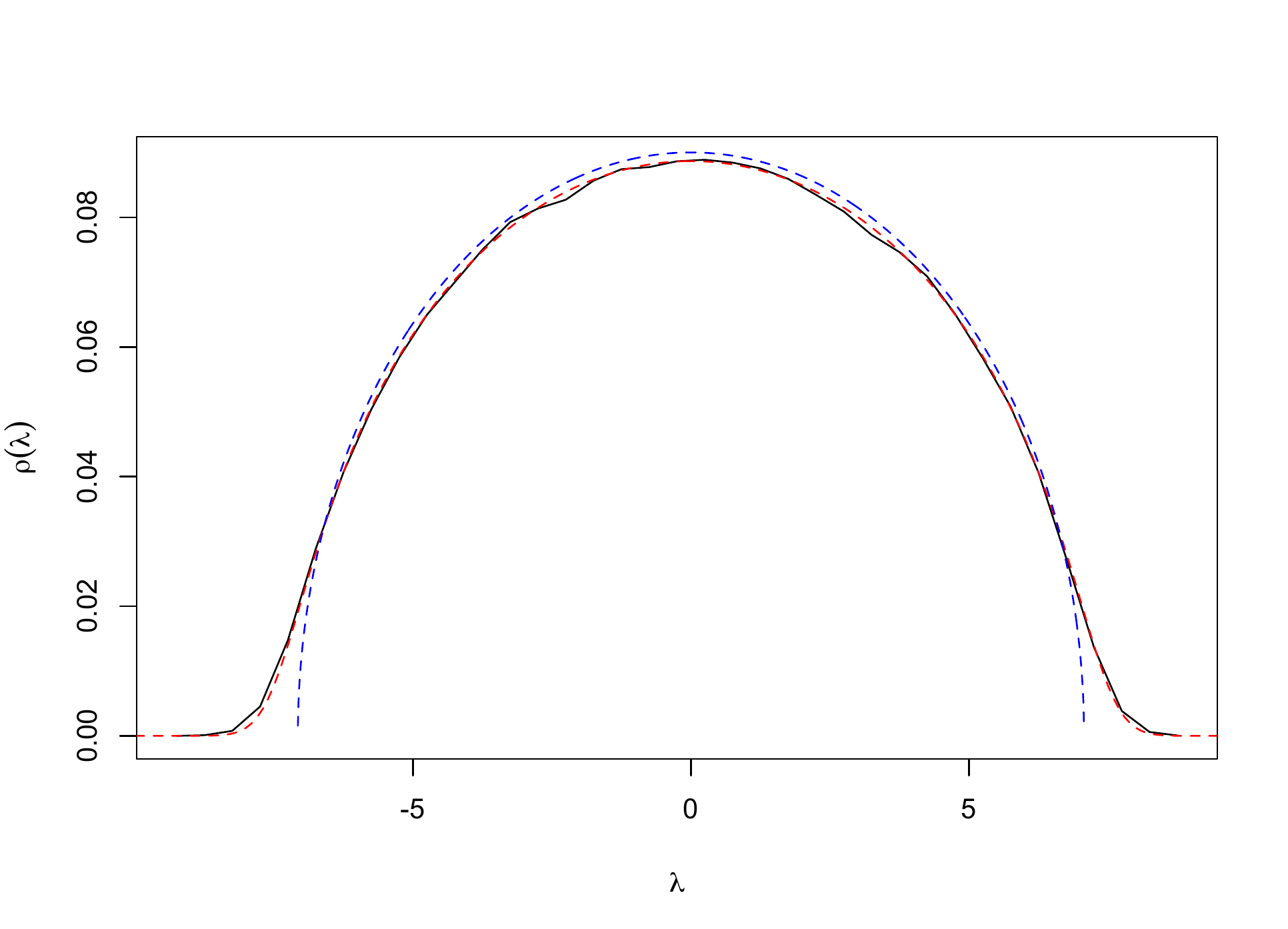}
		\caption{Numerical simulations of the state density of ${\bf M}_n(t=\infty)$ for $p=\beta=1/2,N=50$. Black curve: sample density. 
		Red curve:  {\it Finite size correction} \eqref{wigner_correction}, which coincides almost perfectly with the sample density.
		Blue curve: $N \to \infty$ Wigner semi-circle density. }\label{fig_correction_wigner}
\end{figure}

Let us return to \eqref{RS} for $\beta=\alpha \in (0;2)$. Interestingly, our method allows us to compute the correction to the Wigner semicircle for the spectral density for large but finite $N$ due to the last 
diffusion term, which is usually neglected. Indeed one can solve as above the stationary equation 
of \eqref{RS} keeping every term. This leads to the following {\it corrected spectral density}, valid for large but finite $N$:
\begin{equation}\label{wigner_correction}
\rho(\lambda) = \frac{\sqrt{\alpha}}{\sqrt{2\pi} \Gamma(1+c)} \frac{1}{|D_{-c}(i \sqrt{\alpha} \lambda)|^2}\,,
\end{equation}
where $\alpha= 2/(2-\beta)$ and $c=\beta N/(2-\beta)$.  Numerical verification for this is provided in Fig. \ref{fig_correction_wigner}.

The above discussion can also be formally extended to $-1 \leq c <0$, corresponding to a weakly attracting Coulomb gas (also mentioned in \cite{Vivo}; 
see also \cite{fyodorov} for an application). The system evolves, away from collisions, as in \eqref{dyson_p_hat} 
but with $c\in (-1;0]$ and with infinitesimal {\it ``up-down pushes"} to separate particles when they collide. 
This kind of system has not been rigorously defined in the literature before but it should be possible, at least for a non trivial range of negative values of $c$, 
since the attraction between particles is the same as for Bessel processes of dimension $\delta\in [0;1]$. 
(see \cite[page 3]{Werner}). 
We conjecture that the stationary density for large system 
is again given by the above
Askey-Wimp-Kerov distributions $\rho_c$ but for the parameter range $c\in (-1;0]$. For $c=-1$, the stationary density $\rho_{-1}$ is 
a Dirac mass at $0$. Beyond this level, the attraction is too strong and the gas completely collapses on itself.

%

As a conclusion, we have provided here the first explicit construction of invariant $\beta$-ensembles of random matrices, for arbitrary $\beta \leq 2$. 
The stationary distribution for the eigenvectors is the Haar probability measure on the orthogonal group  if $0< \beta \leq 1$, respectively unitary group if $1<\beta \leq 2$. We have 
found a natural scaling limit that allows one to interpolate smoothly between the Gaussian distribution, relevant for sums of independent random variables, and 
the Wigner semi-circle distribution, relevant for sums of free random matrices. The interpolating limit distributions form a one parameter family that can be explicitly computed.
Let us mention three interesting open problems. First, our alternate Bernoulli process of commuting and free matrix `slices' can probably be done differently, for example by introducing
an Ornstein-Uhlenbeck process on the orthogonal group that mean-reverts towards the identity matrix. By using these matrices $O$ to construct the `slices' as ${\rm d}{\bf M}(t) = O^T \,
{\rm d}\Delta \, O$, one may be able to generate other interesting ensembles by tuning the parameters of the Ornstein-Uhlenbeck process, again interpolating between commuting and free 
addition. Second, it would be interesting to know how the 
eigenvalue spacing distribution smoothly interpolates between the Poisson distribution and Wigner's surmise. Finally, the statistics of the largest 
eigenvalue is also very interesting (and now well known for $\beta>0$, see \cite{Laure,Celine,Celine2,Forrester2}): one should be able to interpolate smoothly, as a function 
of $c$, between the well-known Gumbel distribution of extreme value 
statistics and the Tracy-Widom$(\beta)$ distributions. Whether this can be mapped into a generalized KPZ/Directed polymer problem remains to be seen.

\subsection*{Acknowledgements}
We are grateful to S. N. Majumdar for many important suggestions. 
We thank C. Garban, R. Rhodes, V. Vargas, F. Benaych-Georges, L. Dumaz, R. Chicheportiche, G. Schehr and P. P. Vivo for useful comments and discussions. A. G. and J.-P. B. thank 
the organizers of the NYU-Abu-Dhabi conference in January 2011, where this work was initiated (on way to buying dates).

\end{document}